\documentclass[12pt]{article} \usepackage{amsfonts}\usepackage{amssymb}\usepackage{amsmath}\usepackage[english]{babel}
\begin{document}
\baselineskip=14pt
\pagestyle{plain}
{\Large
\makeatletter
\renewcommand{\@makefnmark}{}
\makeatother
\newcommand{\al}{\alpha_{2n-\theta}}
\newcommand{\be}{\beta_{2n-\theta}}
\newcommand{\alk}{\alpha_{2n_k-\theta}}
\newcommand{\bek}{\beta_{2n_k-\theta}}
\newcommand{\wm}{w(\mu_{n,j},0)w(\mu_{n,j},1)}
\newcommand{\fm}{F_{1,\theta}(\mu_{n,j})}
\newcommand{\fn}{F_{2,\theta}(\mu_{n,j})}
\newcommand{\mun}{\mu_{n,j}}
\newcommand{\on}{o(n^{-1})}
\newcommand{\onm}{o(n^{-m-1})}
\newcommand{\alg}{\varphi_{2n-\theta}}
\newcommand{\beg}{\psi_{2n-\theta}}

\newcommand{\uo}{\stackrel{0}{u}(x)}
\newcommand{\en}{\{e_n\}_{n=1}^\infty}
\newcommand{\dl}{\Delta(\lambda)}

\newcommand{\wnj}{W_{N,j}(\mu)}
\newcommand{\ntm}{|2n-\theta-\mu|}
\newcommand{\pwp}{PW_\pi^-}
\newcommand{\cpm}{\cos\pi\mu}
\newcommand{\spm}{\frac{\sin\pi\mu}{\mu}}
\newcommand{\tet}{(-1)^{\theta+1}}
\newcommand{\dm}{\Delta(\mu)}

\newcommand{\dnm}{\Delta_N(\mu)}
\newcommand{\sni}{\sum_{n=N+1}^\infty}
\newcommand{\muk}{\sqrt{\mu^2+q_0}}
\newcommand{\agt}{\alpha, \gamma, \theta,}
\newcommand{\muq}{\sqrt{\mu^2-q_0}}
\newcommand{\smn}{\sum_{n=1}^\infty}
\newcommand{\lop}{{L_2(0,\pi)}}
\newcommand{\tdm}{\tilde\Delta_+(\mu_n)}
\newcommand{\dmn}{\Delta_+(\mu_n)}
\newcommand{\dd}{D_+(\mu_n)}
\newcommand{\ddd}{\tilde D_+(\mu_n)}
\newcommand{\sn}{\sum_{n=1}^\infty}
\newcommand{\pnn}{\prod_{n=1}^\infty}
\newcommand{\aln}{\alpha_n(\mu)}
\newcommand{\emp}{e^{\pi|Im\mu|}}
\newcommand{\ppp}{\prod_{p=p_0}^\infty}
\newcommand{\ppk}{\sum_{p=p_0}^{\infty}\sum_{k=1}^{[\ln p]}}
\newcommand{\pps}{\sum_{p=p_0}^\infty}
\newcommand{\mln}{m(\lambda_n)}
\newcommand{\mnk}{\mu_{n_k}}
\newcommand{\lnk}{\lambda_{n_k}}
\newcommand{\lqr}{\langle q \rangle}
\newcommand{\mlnk}{m(\lambda_{n_k})}

\newcommand{\sxm}{s(x,\mu)}
\newcommand{\skm}{s(\xi,\mu)}
\newcommand{\cxm}{c(x,\mu)}
\newcommand{\ckm}{c(\xi,\mu)}

\newcommand{\vs}{\vspace{.3cm}}

\centerline {\bf On two-point boundary value problems for }
\centerline {\bf the Sturm-Liouville operator}

\medskip
\medskip
\centerline {\bf Alexander Makin}
\medskip
\medskip

\footnote{2000 Mathematics Subject Classification. 34L10.

Key words and phrases. Sturm-Liouville operator, basis property, root function system.}
\begin{abstract}
In this paper, we study spectral problems for the Sturm-Liouville operator with arbitrary complexvalued potential
$q(x)$ and two-point boundary conditions. All types of mentioned boundary conditions are considered. We ivestigate
in detail the completeness property and the basis property of the root function system.
\end{abstract}

\medskip
\medskip

{\bf 1. Introduction.} The spectral theory of two-point differential operators was
begun by Birkhoff in his two papers [1, 2] of 1908 where he introduced regular boundary conditions
for the first time. It was continued by Tamarkin [3, 4] and Stone [5, 6]. Afterwards
their investigations were developed in many directions. There is an enormous
literature related to the spectral theory outlined above, and we refer to [7-18] and their
extensive reference lists for this activity.

The present communication is a brief survey of results in the spectral theory of
the Sturm-Liouville equation
$$
 u''-q(x)u+\lambda u=0\eqno (1)
$$
with two-point boundary conditions
$$
B_i(u)=a_{i1}u'(0)+a_{i2}u'(\pi)+a_{i3}u(0)+a_{i4}u(\pi)=0, \eqno (2)
$$
where the $B_i(u)$ $(i=1,2)$ are linearly independent forms with arbitrary
complex-valued coefficients and $q(x)$ is an arbitrary complex-valued
function of class $L_1(0,\pi)$.

Our main focus is on the non-self-adjoint case, and, in particular, the case when boundary conditions are degenerate.  We will study the completeness
property and the basis property of the root function system of operator (1), (2).
The convergence of spectral expansions is investigated only in classical sense, i.e.
the question about the summability of divergent series by a generalized method is not considered.

 We consider the operator $Lu=u''-q(x)u$ as a linear operator on $\lop$
with the domain $D(L)=\{u\in\lop\, |\, u(x), u'(x)$ are absolutely
continuous on $[0,\pi]$, $u''-q(u)u\in\lop$, $B_i(u)=0$ $(i=1,2)\}$.

By {\it an eigenfunction} of the operator $L$ corresponding to an eigenvalue
 $\lambda\in\mathbb{C}$ we mean any function $\uo\in D(L)$ $(\uo\not\equiv0)$
which satisfies the equation
$$
L\stackrel{0}{u}+\lambda\stackrel{0}{u}=0
$$
almost everywhere on $[0,\pi]$.

By {\it an associated function} of the operator $L$ of order $p$ $(p=1,2, \ldots)$
corresponding to the same eigenvalue $\lambda$ and the eigenfunction $\uo$
we mean any function $\stackrel{p}{u}(x)\in D(L)$ which satisfies the equation
$$
L\stackrel{p}{u}+\lambda \stackrel{p}{u}=\stackrel{p-1}{u}
$$
almost everywhere on $[0,\pi]$. One can also say that an eigenfunction $\uo$ is
an associated function of zero order. The set of all eigen- and associated functions
(or root functions) corresponding to the same eigenvalue $\lambda$ together with
the function $u(x)\equiv0$ forms a root linear manifold. This manifold is called
a root subspace if its dimension is finite.

Let the set of the eigenvalues of the operator $L$ be countable and
all root linear manifolds be root subspaces. Let us choose a basis in each root
subspace. Any system $\{u_n(x)\}$ obtained as the union of chosen bases of all the root
subspaces is called {\it a system of eigen- and associated functions} (or root function
system) of the operator $L$.

The main purpose of this lecture is to study the basis property of the root function
system of the operator $L$. Before starting our investigation we must verify
completeness of the root function system in $\lop$.

It is convenient to write conditions (2) in the matrix form
$$
A=\left(
\begin{array}{cccc}
a_{11}&a_{12}&a_{13}&a_{14}\\
a_{21}&a_{22}&a_{23}&a_{24}
\end{array}
\right)
$$
and denote the matrix composed of the ith and jth columns of $A$
  $(1\le i<j\le4)$
by $A(ij)$; we set $A_{ij}=det A(ij)$.

Denote by $c(x,\mu), s(x,\mu)$ $(\lambda=\mu^2)$ the fundamental system of solutions to equation (1)
with the initial conditions $c(0,\mu)=s'(0,\mu)=1$, $c'(0,\mu)=s(0,\mu)=0$.
The eigenvalues of problem (1), (2) are the roots of the characteristic determinant
$$
\Delta(\mu)
=
\left|
\begin{array}{cc}
B_1(c(x,\mu))&B_1(s(x,\mu))\\
B_2(c(x,\mu))&B_2(s(x,\mu))
\end{array}
\right|.
$$
Simple calculations show that
$$
\Delta(\mu)=-A_{13}-A_{24}+A_{34}s(\pi,\mu)-A_{23}s'(\pi,\mu)-A_{14}c(\pi,\mu)-A_{12}c'(\pi,\mu).
$$
It is easily seen that if $q(x)\equiv0$ then the characteristic determinant $\Delta_0(\mu)$
of the corresponding problem (1), (2) has the form
$$
 \Delta_0(\mu)=-A_{13}-A_{24}+A_{34}\frac{\sin\pi\mu}{\mu}-(A_{23}+A_{14})\cos\pi\mu
 +A_{12}\mu\sin\pi\mu.
 $$

Boundary conditions (2) are called {\it nondegenerate} if they satisfy one of the
following relations:
$$
1) A_{12}\ne0,\quad2) A_{12}=0, A_{14}+A_{23}\ne0,\quad3) A_{12}=0,
 A_{14}+A_{23}=0, A_{34}\ne0.
$$
Evidently, boundary conditions (2) are nondegenerate iff $\Delta_0(\mu)\ne const$.

Notice, that for any
nondegenerate boundary conditions an asymptotic representation for the characteristic determinant $\Delta(\mu)$ as $|\mu|\to\infty$ one can find in [10].

{\bf Theorem 1} ([10]). {\it For any
 nondegenerate conditions the
spectrum of problem (1), (2) consists of a countable set $\{\lambda_n\}$
of eigenvalues with only one limit point $\infty$, and the dimensions of the
corresponding root subspaces are bounded by one constant. The system $\{u_n(x)\}$
of eigen- and associated functions is complete and minimal in $L_2(0,\pi)$; hence,
it has a biorthogonally dual system $\{v_n(x)\}$.}

 For convenience,
we introduce numbers $\mu_n$, where $\mu_n$ is the square root of $\lambda_n$
with nonnegative real part.

It is known that nondegenerate conditions can be divided into three classes:

1) strengthened regular conditions;

2) regular but not strengthened regular conditions;

3) irregular conditions.

The definitions are given in [8]. These three cases should be considered
separately.

{\bf 2. Strengthened regular conditions.} Let boundary conditions (2) belong to
class 1). According to
[8], this is equivalent to the fulfillment one of the following conditions:
$$
\begin{array}{c}
A_{12}\ne0;\quad A_{12}=0,
 A_{14}+A_{23}\ne0,
 A_{14}+A_{23}\ne\mp(A_{13}+A_{24});\\
 A_{12}=0, A_{14}+A_{23}=0,
A_{13}+A_{24}=0, A_{13}=A_{24},
 A_{34}\ne0.
\end{array}
$$
It is well known that,
 all but finitely many eigenvalues $\lambda_n$ are simple
(in other words, they are asymptotically simple),
and the number of associated functions is finite. Moreover, the $\lambda_n$
are separated in the sense that there exists a constant $c_0>0$ such that, for
any sufficiently large different numbers $\lambda_k$ and $\lambda_m$, we have
$$
|\mu_k-\mu_m|\ge c_0. \eqno (3)
$$

{\bf Theorem 2.}{\it
The system of root functions $\{u_n(x)\}$ forms a Riesz basis in $L_2(0,\pi)$.}

This statement was proved in [21], [22] and [9, Chapter XIX].

Class 1) contains many types of boundary conditions, for example, the Dirichlet
boundary conditions $u(0)=u(\pi)=0$, the Newmann boundary conditions $u'(0)=u'(\pi)=0$,
the Dirichlet-Newmann boundary conditions $u(0)=u'(\pi)=0$ and others.

{\bf 3. Regular but not strengthened regular conditions.}
Let boundary conditions belong to class 2).
According to [8], this is equivalent to the fulfillment of the conditions
$$
A_{12}=0, \quad A_{14}+A_{23}\ne0,\quad A_{14}+A_{23}=(-1)^{\theta+1}(A_{13}+A_{24}),\eqno(4)
$$
where $\theta=0,1$. It is well known [10] that the eigenvalues of problem (1), (2)
form two series:
$$
\lambda_0=\mu_0^2, \quad \lambda_{n,j}=(2n+o(1))^2 \eqno(5)
$$
(if $\theta=0$) and
$$
\lambda_{n,j}=(2n-1+o(1))^2 \eqno(6)
$$
(if $\theta=1$). Here, in both cases, $j=1,2$ and $n=1,2,\ldots$. We denote
$\mu_{n,j}=\sqrt{\lambda_{n,j}}=2n-\theta+o(1)$. It follows from [8] that
asymptotic formulas (6) and (7) can be refined. Specifically,
$$
\mu_{n,j}=2n-\theta+O(n^{-1/2}).
$$
Obviously, $|\mu_{n,1}-\mu_{n,2}|=O(n^{-1/2})$; i.e. $\mu_{n,1}$ and $\mu_{n,2}$
become infinitely close to each other as $n\to\infty$. If $\mu_{n,1}=\mu_{n,2}$
for all $n$, except, possibly, a finite set, then the spectrum of problem (1), (2)
is called {\it asymptotically multiple}. If the set of multiple eigenvalues is
finite, then the spectrum of problem (1), (2) is called {\it asymptotically simple}.

There exist numerous examples when
the number of multiple eigenvalues is
finite or infinite, and the total number of associated functions is finite
or infinite also. We see that separation condition
(3) never holds. Depending on the particular form of the boundary
conditions and the potential $q(x)$
the system of root functions may have or may not have the basis property [17], [22], [23], and
even for fixed boundary conditions, this property may appear or disappear under
arbitrary small variations of the coefficient $q(x)$ in the corresponding metric [24].
Thus, the considered case is much more complicated than the previous one, so we will
study it in detail.

For any problem (1), (2) let $Q$ denote the set of potentials $q(x)$ from the
class $L_1(0,\pi)$ such that the system of root functions forms a Riesz basis
in $L_2(0,\pi)$, $\bar Q=L_1(0,\pi)\setminus Q$.

To analyze this class of problems, it is reasonable [12] to divide
conditions (2) satisfying (4) into three types:

I) $A_{14}=A_{23}$, $A_{34}=0$;

II) $A_{14}=A_{23}$, $A_{34}\ne 0$;

III) $A_{14}\ne A_{23}$

The eigenvalue problem for operator (1) with boundary conditions of type
I, II, or III, is called the problem of type I, II, or III, respectively.

At first we consider the problems of type I.
It was shown in [12] that any boundary conditions of type I are equivalent
to the boundary conditions specified by the matrix
$$
A=\left(
\begin{array}{cccc}
1&(-1)^{\theta+1}&0&0\\
0&0&1&(-1)^{\theta+1}
\end{array}
\right),
$$
i.e., to periodic or antiperiodic boundary conditions. These boundary conditions
are selfadjoint.

{\bf Theorem 3} ([25]). {\it The sets $Q$ and $\bar Q$ are everywhere dense in $L_1(0,\pi)$.}

Recently, (see [26-37] and their extensive reference lists)
by a number of authors, a very nice theory of the problems of type I was built.

Let us consider the problems of type II. It was also established in [12]
that any boundary conditions
of type II are equivalent to the boundary conditions specified by the matrix

$$
A=\left(
\begin{array}{cccc}
1&-1&0&a_{14}\\
0&0&1&-1
\end{array}
\right) \quad \mbox{or}\quad
A=\left(
\begin{array}{cccc}
1&1&0&a_{14}\\
0&0&1&1
\end{array}
\right),
$$
where $a_{14}\ne0$ in both cases. If $a_{14}$ is a real number and
$q(x)$ is a real function, then the corresponding boundary value problem is selfadjoint.

{\bf Theorem 4} ([38]). {\it If $A_{14}=A_{23}$ and $A_{34}\ne0$, then the system $\{u_n(x)\}$
forms a Riesz basis in $L_2(0,\pi)$, and the spectrum is asymptotically simple.}
Denote by $\{v_n(x)\}$ the biorthogonally dual system.
The key point in the proof of Theorem 4 is obtaining the estimate
$$
\max_{(x,\xi)\in[0,\pi]\times[0,\pi]}|u_n(x)\overline{v_n(\xi)}|\le C,\eqno(7)
$$
which is valid for any number $n$.
It follows from (7) and [39] that the system $\{u_n(x)\}$ forms a Riesz basis in
$L_2(0,\pi)$.

A comprehensive description of boundary conditions of type III was
given in [12].
In particular, it is known that all of them are non-self-adjoint.

{\bf Theorem 5} ([38]). {\it If
$
A_{14}\ne A_{23}
$, then the system of root functions
$\{u_n(x)\}$ of problem (1), (2) is a Riesz basis in $L_2(0,\pi)$
if and only if the spectrum is asymptotically multiple
.}

Thus, we have established that for problems of type III the question
about the basis property for the system of eigen- and associated functions is
reduced to the question about asymptotic multiplicity of the spectrum.
The presence of this property depends essentially on the particular form of the
boundary conditions and the function $q(x)$.

{\bf Theorem 6} ([40, 41]). {\it If $A_{14}\ne A_{23}$, then, for any function $q(x)\in L_2(0,\pi)$
and any $\varepsilon>0$,
there exists a function
$\tilde q(x)\in L_2(0,\pi)$ such that $||q(x)-\tilde q(x)||_{L_2(0,\pi)}<\varepsilon$
and
problem (1), (2) with the potential $\tilde q(x)$ has an asymptotically
multiple spectrum.}

For $A_{14}=A_{23}$ and $A_{34}=0$, a similar proposition was deduced in [42].

Theorems 3, 4, 6 and the results of
[43] imply that the whole class of regular but not strengthened regular
boundary conditions splits into two subclasses (a) and (b). Subclass (a)
coincides with the second type of boundary conditions and is characterized by the fact
that the system of root functions of problem (1), (2) with boundary conditions
from this subclass forms a Riesz basis in $L_2(0,\pi)$ for any potential $q(x)\in L_1(0,\pi)$;
i.e. $Q=L_1(0,\pi)$, $\bar Q=\emptyset$. We will see below that
boundary conditions from the subclass (a) are the only boundary conditions
(in addition to strengthened regular ones) that ensure the Riesz basis property of
the system of root functions for any potential $q(x)\in L_1(0,\pi)$.

Subclass (b) contains the remaining regular but not strengthened regular
boundary conditions. An entirely different situation takes place in this case.
For any problem with boundary conditions from this subclass, the sets
$Q$ and $\bar Q$ are dense everywhere
in $L_1(0,\pi)$.

{\bf 4. Irregular conditions.} Let boundary conditions (2) belong to
class 3). According to
[8, 12], this is equivalent to the fulfillment one of the following conditions:
$$
\begin{array}{c}
A_{12}=0,\quad
 A_{14}+A_{23}=0,\quad
 A_{13}+A_{24}=0,\quad
 A_{13}\ne A_{24},\quad
 A_{34}\ne0;
\\
 A_{12}=0,\quad A_{14}+A_{23}=0,\quad
A_{13}+A_{24}\ne0,\quad
 A_{34}\ne0.
\end{array}
$$

According to [12], any boundary
conditions of the considered class
are equivalent to the boundary conditions determined by the matrix
$$
A=\left(
\begin{array}{cccc}
1&\pm1&0&b_0\\
0&0&1&\mp1
\end{array}
\right),\quad
\mbox{where}\quad b_0\ne0,
$$
or

$$
A=\left(
\begin{array}{cccc}
1&b_1&0&b_0\\
0&0&1&-b_1
\end{array}
\right), \quad \mbox{where} \quad b_1\ne\pm1,\quad b_0\ne0,
$$
or
$$
A=\left(
\begin{array}{cccc}
0&1&a_0&0\\
0&0&0&1
\end{array}
\right),\quad \mbox{where} \quad a_0\ne0.
$$

In case 3), as well as in case 1),
all but finitely many eigenvalues $\lambda_n$ are simple, the number of
associated functions is finite, and separation condition
(3) holds. However, the system $\{u_n(x)\}$ never forms even a usual basis in $L_2(0,\pi)$,
because $||u_n||_{L_2(0,\pi)}||v_n||_{L_2(0,\pi)}\to\infty$ as $n\to\infty$. Here $\{v_n(x)\}$ is the biorthogonally dual system.This
case was investigated in [5], [6], [44].

{\bf 5. Degenerate conditions.} Let boundary conditions (2) be degenerate. According to
[10, 12], this is equivalent to the fulfillment of the following conditions:

$$
A_{12}=0, \quad A_{14}+A_{23}=0,\quad A_{34}=0.
$$

According to [12], any boundary
conditions of the considered class
are equivalent to the boundary conditions determined by the matrix
$$
A=\left(
\begin{array}{cccc}
1&d&0&0\\
0&0&1&-d
\end{array}
\right),\quad
\mbox{or}\quad	A=\left(
\begin{array}{cccc}
0&1&0&0\\
0&0&0&1
\end{array}
\right).
$$
If in the first case $d=0$ then for any potential $q(x)$ we have the
initial value problem (the Cauchy problem) which has no eigenvalues. The same
situation takes place in the second case.

Further we will consider the first case if $d\ne0$. Then the boundary conditions can be written
in more visual form
$$
u'(0)+du'(\pi)=0,\quad u(0)-du(\pi)=0. \eqno(8)
$$
{\bf 5.1. Spectrum.}

By $PW_\sigma$ we denote the class of entire functions $f(z)$ of exponential type
$\le\sigma$ such that $||f(z)||_{L_2(R)}<\infty$, and by $PW_\sigma^-$
we denote the set of odd functions in $PW_\sigma$.

By performing simple manipulations, we obtain the relation
$$
\begin{array}{c}
\Delta(\mu)=\frac{d^2-1}{d}+c(\pi,\mu)-s'(\pi,\mu)=
\frac{d^2-1}{d}+\int_0^\pi K(t)\frac{\sin\mu t}{\mu}dt=\\
=\frac{d^2-1}{d}+\frac{f(\mu)}{\mu},
\end{array}
$$
where $K(t)\in L_1(0,\pi)$. If $q(x)\in L_2(0,\pi)$ then $K(t)\in L_2(0,\pi)$ and $f(\mu)\in PW^-_\pi$.
Notice, that simple calculations show that if $d=\pm1$ and
$q(x)\equiv0$ then any $
 \lambda\in\mathbb{C}$ is an eigenvalue of infinite multiplicity. This abnormal
example constructed by Stone illustrates the difficulty of investigation of problems with boundary conditions of the
considered class.

It is well known that the characteristic determinant $\Delta(\mu)$ of problem (1), (8)
is an entire function of the parameter $\mu$, consequently, for the operator (1), (8) we have only the following possibilities:

1) the spectrum is absent;

2) the spectrum is a finite nonempty set;

3) the spectrum is a countable set without finite limit points;

4) the spectrum  fills the entire complex plane.

 One can prove [45] that case 2) is impossible.
It is known that $c(\pi,\mu)-s'(\pi,\mu)\equiv0$ if and only if the function
$$
Q(x)=q(x)-q(\pi-x)=0 \eqno(9)
$$
almost everywhere  on the segment $[0,\pi]$. Evidently, $\lim_{\mu\to\infty}(c(\pi,\mu)-s'(\pi,\mu))=0$.
 Hence it follows that if $c(\pi,\mu)-s'(\pi,\mu)\equiv C$ then $C=0$. We see that case 1) takes place if and only if condition (9) holds and $d\ne\pm1$, and case 4) takes place if and only if condition (9) holds and $d=\pm1$. If condition (9) does not hold we have case 3).

 {\bf 5.2. Completeness.}

  Completeness of the root function system of problem (1), (8)
 was investigated in [46-47]. The main result of the mentioned papers is:

{\bf Theorem 7} ([47]). {\it If $q(x)\in C^k[0,\pi]$ for some $k=0,1\ldots$ and $q^{(k)}(0)\ne(-1)^kq^{(k)}(\pi)$,
then the system of root functions is complete in the space $L_p(0,\pi)$ if $1\le p<\infty$.}

It follows from [47]
that
depending on the potential $q(x)$
the system of root functions may have or may not have the completeness property,
moreover, this property may appear or disappear under
arbitrary small variations of the coefficient $q(x)$ in the corresponding
metric even for fixed boundary conditions.

If the conditions $d^2\ne1$ and $q(x)\in C[0,\pi]$ hold necessary and sufficient conditions of the completeness of root function system of problem
(1), (8) were found in
[48].

{\bf Theorem 8}  ([49]). {\it If for a number $\rho>0$
$$
\lim_{h\to0}\frac{\int_{\pi-h}^{\pi}Q(x)dx}{h^\rho}=\nu,
$$
and $\nu\ne0$,
then the root function system of problem (1), (8) is complete in the space $L_p(0,\pi)$ if $1\le p<\infty$.

}

Since for a wide class of potentials $q(x)$ the root function system of problem (1), (8) is complete in $L_2(0,\pi)$ one can set a question whether the mentioned system forms a basis.

 {\bf 5.3. Basis property.}

Recently, it was proved in [50] that the root function system never forms an unconditional
basis in $\lop$ if multiplicities of the eigenvalues are uniformly bounded by some constant.
Moreover, under the condition mentioned above it was established there that if the eigen- and associated function system of general ordinary differential operator with two-point boundary conditions  forms an unconditional basis then the
boundary conditions  are regular. Article [50]
was published in 2006. At that time it was unknown whether there exists a potential $q(x)$ providing unbounded
growth of multiplicities of the eigenvalues.

{\bf 5.4. Inverse problem.}

 However, in 2010 in [45] an example of a potential $q(x)$ for which
the characteristic determinant has the roots of arbitrary high multiplicity was constructed. Hence, the corresponding root
function system $\{u_n(x)\}$ contains associated functions of arbitrary high order. It means, that paper [50] does not give the definitive solution of basis property problem.
Below we will show a method to construct a potential $q(x)$ providing unbounded growth of multiplicities of eigenvalues.

{\bf Theorem 9} ([51]). {\it Suppose that a function $v(\mu)$ can be represented in the form
$$
v(\mu)=\gamma+\frac{f(\mu)}{\mu},\eqno(10)
$$
where $\gamma$ is some complex number, the function$f(\mu)\in PW_\pi^-$
satisfies the condition
$$
\int_{-\infty}^{\infty}|\mu^{m}f(\mu)|^2d\mu<\infty,
$$ where
$m$ is a nonnegative integer number.
Then there exists a function $q(x)\in W_2^m(0,\pi)$ such that the characteristic determinant $\Delta(\mu)$ of problem (1), (8), where either
 $d=(\gamma+\sqrt{\gamma^2+4})/2$ or
 $d=(\gamma-\sqrt{\gamma^2+4})/2$
and with the potential  $q(x)$ is identically equal to the function $v(\mu)$.}

Therefore, Theorem 9 reduces the problem on the structure of the spectrum of problem (1), (8) with degenerate boundary conditions to the problem on the expansion of a function of the form (10) into a canonical product.

{\bf 5.5. Nontrivial examples.}

{\it Example 1.} Let us define a sequence $\{a_k\}$ $(k=1,2,\ldots)$ in this way: $a_1=1, a_2=3, a_3=5$, $a_{k+1}=a_k+2p$,
if $2^p<k<2^{p+1}$ $(p=1,2,\ldots)$, and $a_{k+1}=a_k+(a_k-a_{k-1})+2$, if $k=2^p$ $(p=2,3,\ldots)$. Set
$$
F_1(\mu)=\prod_{k=1}^\infty\left(1-\frac{\mu^2}{a_k^2}\right)^{a_{k+1}-a_k-\delta_k},
$$
where $\delta_k=0$, if $k\ne2^p$, and $\delta_k=1$, if $k=2^p$, $(p=2,3,\ldots)$.

{\bf Theorem 10} ([52]). {\it For any real $x$  the following inequality holds
$$
|F_1(x)|\le C_3(|x|+1)^{M_1},
$$
where $M_1$ is a sufficiently large number.}

 Denote
$$
f_1(\mu)=\mu\prod_{k=M_1+1}^\infty\left(1-\frac{\mu^2}{a_k^2}\right)^{a_{k+1}-a_k-\delta_k}.
$$
It is easily shown that $f_1(\mu)\in PW_\pi^-$. This, together with Theorem 9 implies that there exists a potential $q_1(x)\in L_2(0,\pi)$, such that for the characteristic determinant $\Delta_1(\mu)$ of problem
$$
u''-q_1(x)u+\lambda u=0,\quad u'(0)+ du'(\pi)=0,\quad
u(0)-du(\pi)=0
 \eqno (11)
 $$
$(d=\pm1)$ we have the equality
$$
\Delta_1(\mu)=f_1(\mu)/\mu.
$$
It follows from the definition of sequence $\{a_k\}$ that multiplicities of zeros $a_k$ of constructed above function $f_1(\mu)$ monotonically not decrease and tend to infinity as $k\to\infty$. Therefore, the eigenvalues $\lambda_n=\mu_n^2$ of problem (11) have the desired property: their multiplicities  $m(\lambda_n)$ tend to infinity and the corresponding root function system contains associated functions of arbitrary high order, i.e. the dimensions of root subspaces infinitely grow. Moreover, the following inequality takes place
$$
c_1\ln\mu_n\le m(\lambda_n)\le  c_2\ln\mu_n.
$$

{\bf Theorem 11} ([52]). {\it
 The root function system $\{u_n(x)\}$ of problem (11) is complete in $L_2(0,\pi)$.}

{\it Example 2.}
Denote $\tilde a_k=a_k-\alpha_k+i\beta_k$, $(k=1,2,\ldots)$ where $\alpha_k=a_k-\sqrt{a_k^2-\beta_k^2}$, $\beta_k=(a_k-a_{k-1})/10)$ $(a_0=0)$.
Denote $h_k=a_{k+1}-a_k-\delta_k$,
$$
F_2(\mu)=\prod_{k=1}^\infty\left(1-\frac{\mu^2}{\tilde a_k^2}\right)^{[h_k/2]}
\left(1-\frac{\mu^2}{\bar{\tilde a}_k^2}\right)^{[h_k/2]}\left(1-\frac{\mu^2}{\tilde a_k^2}\right)^{h_k-2[h_k/2]}.
$$
{\bf Theorem 12} ([53]). {\it For any real $x$  the following inequality holds
$$
|F_2(x)|\le C_3(|x|+1)^{M_2},
$$
where $M_2$ is a sufficiently large number.}

 Denote
$$
f_2(\mu)=\mu\prod_{k=M_2+1}^\infty\left(1-\frac{\mu^2}{\tilde a_k^2}\right)^{[h_k/2]}
\left(1-\frac{\mu^2}{\bar{\tilde a}_k^2}\right)^{[h_k/2]}\left(1-\frac{\mu^2}{\tilde a_k^2}\right)^{h_k-2[h_k/2]}.
$$

It is easily shown that $f_2(\mu)\in PW_\pi^-$. This, together with Theorem 9 implies that there exists a potential $q_2(x)\in L_2(0,\pi)$, such that for the characteristic determinant $\Delta_2(\mu)$ of problem
$$
u''-q_2(x)u+\lambda u=0,\quad u'(0)+ du'(\pi)=0,\quad
u(0)-du(\pi)=0
 \eqno (12)
 $$
$(d=\pm1)$ we have the equality
$$
\Delta_2(\mu)=f_2(\mu)/\mu.
$$
It follows from the definition of sequence $\{\tilde a_k\}$ that multiplicities of zeros $\tilde a_k, \bar{\tilde a}_k$ of constructed above function $ f_2(\mu)$ monotonically not decrease and tend to infinity as $k\to\infty$. Therefore, the eigenvalues $\tilde\lambda_n=\tilde\mu_n^2$ of problem (12) have two properties: their multiplicities $m(\tilde\lambda_n)$ tend to infinity, hence, the corresponding root function system contains associated functions of arbitrary high order, and $|Im\tilde\mu_n|\to\infty$ as $n\to\infty$. Moreover, the following two inequalities hold
$$
c_1\ln|\tilde\mu_n|\le m(\lambda_n)\le  c_2\ln|\tilde\mu_n|,  \quad c_1|Im \tilde\mu_n|\le m(\tilde\lambda_n)\le  c_2\ln|\tilde\mu_n|.
$$

{\bf Theorem 13} ([53]). {\it
 The root function system $\{u_n(x)\}$ of problem (12) is complete in $L_2(0,\pi)$.}

For any problem (1), (8) let $\Omega$ denote the set of potentials $q(x)$ from the
class $L_1(0,1)$ such that the system of root functions contains associated functions of arbitrary high order, $\bar \Omega=L_1(0,\pi)\setminus \Omega$.

{\bf Theorem 14} ([54]). {\it The sets $\Omega$ and $\bar\Omega$ are everywhere dense in $L_1(0,\pi)$.}

Since for a wide class of potentials $q(x)$ the root function system of problem (1), (8) is complete in $L_2(0,\pi)$ one can set a question whether the mentioned system forms a basis.

{\bf 5.6. Again on the basis property.}

Let $\lambda_n=\mu_n^2$ $(Re\mu_n\ge0, n=1,2,\ldots)$ be the eigenvalues of problem (1), (8) numbered neglecting their multiplicities in nondecreasing order of absolute value. By $\mln$ we denote the multiplicity of an eigenvalue $\lambda_n$.

{\bf Theorem 15} ([55-57]). {\it If

$$
\lim_{n\to\infty}\frac{\mln}{\sqrt{|\mu_n|}}=0,\eqno(13)
$$
then the system of eigenfunctions and associated functions of problem (1), (8) is not a basis in $L_2(0,\pi)$.}

Clearly, since Theorem 15 contains supplementary condition (13), it does not give the definitive solution of the basis property problem.
If this condition does not hold then the mentioned problem has not been solved.
\medskip

\centerline {\bf Acknowledgement}
\medskip

This work was supported by the Russian Foundation for Basic Research, project No. 13-01-00241.
\medskip
\medskip
\medskip

\centerline {\bf REFERENCES}
\medskip
\medskip
\medskip

[1] G.D. Birkhoff. On the asymptotic character of the solutions of certain linear
differential equations containing a parameter. Trans. Amer. Math. Soc. {\bf 9} (1908), 219-231.

[2] G.D. Birkhoff. Boundary value and expansions problems of ordinary linear
differential equations. Trans. Amer. Math. Soc. {\bf 9} (1908), 373-395.

[3] J. Tamarkin. Sur Quelques Points de la Theorie des Equations Differentielles
Lineaires Ordinaires et sur la Generalisation de la serie de Fourier. Rend. Circ.
Matem. Palermo {\bf 34} (1912), 345-382.

[4] J. Tamarkin. Some general problem of the theory of ordinary linear differential
equations and expansions of an arbitrary function in series of fundamental
functions. Math. Z. (1927), 1-54.

[5] M.H. Stone. A comparison of the series of Fourier and Birkhoff.
Trans. Amer. Math. Soc. {\bf 29} (1926), 695-761.

[6] M.H. Stone. Irregular differential systems of order two and the related expansions
problems. Trans. Amer. Soc. {\bf 29} (1927), 23-53.

[7] E.A. Coddington and N. Levinson. Theory of Ordinary Differential Equations.
McGraw-Hill, New-York, 1955.

[8] M.A. Naimark. Linear Differential Operators (in Russian). Nauka, Moscow, 1969; English transl.:
 Ungar, New-York, 1967.

[9] N. Danford and J.T. Schwartz. Linear Operators, Part III, Spectral Operators, Wiley,
New-York, 1971.

[10] V.A. Marchenko. Sturm-Liouville Operators and Their Applications. Kiev, 1977
(in Russian); English transl.: Birkh\"{a}user, Basel, 1986.

[11] P. Lang and J. Locker. Spectral theory of two-point differential operators
determined by $-D^2$. I. Spectral properties. J. Math. Anal. Appl. {\bf 141} (1989), 538-558.

[12] P. Lang and J. Locker. Spectral theory of two-point differential operators
determined by $-D^2$. II. Analysis of cases. J. Math. Anal. Appl. {\bf 146} (1990), 148-191.

[13]  J. Locker. The spectral theory of second order two-point differential
operators: I. A priori estimates for the eigenvalues and completeness.
Proc. Roy. Soc. Edinburgh, Sect. A. {\bf 121} (1992), 279-301.

[14]  J. Locker. The spectral theory of second order two-point differential
operators: II. Asymptotic expansions and the characteristic determinant.
J. Differential Equations. {\bf 114} (1994), 272-287.

[15]  J. Locker. The spectral theory of second order two-point differential
operators: III. The eigenvalues and their asymptotic formulas.
Rocky Mountain J. Math. {\bf 26} (1996), 679-706.

[16]  J. Locker. The spectral theory of second order two-point differential
operators: IV. The associated projections and the subspace $S_\infty(L)$.
Rocky Mountain J. Math. {\bf 26} (1996), 1473-1498.

[17] V.A. Il'in and L.V. Kritskov. Properties of spectral expansions corresponding
to non-self-adjoint differential
operators. J. of Mathematical Sciences. {\bf 116} (2003), 3489-3550.

[18] J. Locker. Spectral Theory of Non-self-adjoint Two-point differential
operators. Math. Surveys Monogr., vol. 192, North-Holland, Amsterdam, 2003.

[19] S.G. Krein. Functional Analysis (in Russian). Nauka, Moscow, 1972.

[20] I. Ts. Gokhberg and M.G. Krein. Introduction to the Theory of Linear Not Self-Adjoint
Operators (in Russian). Nauka, Moscow, 1965.

[21] V.P. Mikhailov. On Riesz bases in $L^2(0,1)$. Dokl. Akad. Nauk SSSR {\bf 144}
(1962), 981-984 (in Russian).

[22] G.M. Kesel'man. On the unconditional convergence of eigenfunction expansions
of certain differential operators. Izv. Vyssh. Uchebn. Zaved. Mat. {\bf 39} (1964),
82-93 (in Russian).

[23] P.W. Walker. A nonspectral Birkhoff-regilar differential operator.
Proc. Amer. Math. Soc. {\bf 66} (1977), 187-188.

[24] V.A. Il'in. On a connection between the form of the boundary conditions
and the basis property and the equiconvergence with a
trigonometric series in root functions of a nonselfadjoint operator. Differ. Uravn. {\bf 30} (1994),
 1516-1529 (in Russian).

[25] A.S. Makin.  Convergence of Expansions in the Root Functions of Periodic
	Boundary Value Problems. Doklady Mathematics, 73, No. 1, 71-76 (2006);
	translation from Doklady Acad. Nauk, 406, No. 4, 452-457 (2006) (Russian).

[26] F. Gesztesy, V. Tkachenko. When is a non-self-adjoint Hill operator a spectral operator of scalar type? C. R. Math. Acad. Sci. Paris, {\bf 343} (2006), 239-242.

[27] P. Djakov, B. Mytyagin. Instability zones of one-dimensional periodic Schro\"{o}dinger  and Dirac operators. (Russian) Uspekhi
Math. Nauk {\bf 61}, (2006), 77-182; translation in Russian Math. Surveys, {\bf 61}, (2006), 663-766.

[28] A.A. Shkalikov  and O.A. Veliev.  On the Riesz basis property of eigen- and associated functions of periodic and antiperiodic Sturm-Liouville problems. Math. Notes, {\bf 85} (2009), 647-660.

[29] F. Gesztesy, V. Tkachenko.  A criterion for Hill operators to be spectral operators of scalar type, J. Analyse Math. {\bf 107} (2009), pp. 287-353.

[30] P. Djakov, B. Mytyagin. Convergence of spectral decompositions of Hill operators with trigonometric potentials. Math. Ann.
{\bf 351}, (2011), 509-540.

[31] P. Djakov and B. Mityagin. Convergence of spectral decompositions of Hill operators with trigonometric polynomial potentials, Math. Ann., {\bf 351} (3) (2011), pp. 509-540.

[32] P. Djakov, B. Mytyagin. Convergence of spectral decompositions of Hill operators with trigonometric polynomials as potentials.
Doklady Mathematics, {\bf 436}, (2011), 11-13 ; translation from Doklady Akad. Nauk,  {\bf 83}, (2011), 5-7
(Russian).

[33] F. Gesztesy, V. Tkachenko.  A Schauder and Riesz basis criterion for non-self-adjoint Schr\"{o}dinger operator with periodic and anti-periodic boundary conditions. Journal of Differential Equations, {\bf 253} (2012), 400-437.

[34] P. Djakov, B. Mytyagin. Criteria for existence of Riesz bases consisting of root functions of Hill and 1D Dirac operators,
 J. Funct. Anal. {\bf 263}, (2012), 2300-2332.

[35] P. Djakov, B. Mytyagin. Equiconvergence of spectral decompositions of Hill operator, (Russian), Dokl. Akad. Nauk  {\bf 445}, (2012)
 498-500; translation in Dokl. Math. {\bf 86} (2012) 542-544.

[36] P. Djakov, B. Mytyagin.  Equiconvergence of spectral decompositions of Hill-Schr/"{o}dinger operators. J. Differential Equations,
{\bf 255}, (2013), 3233-3283.

[37] P. Djakov, B. Mytyagin. Riesz basis property of Hill operators with potentials in weighted spaces. Trans. Moscow Math. Soc.
2014, 151-172.

[38] A.S. Makin.  On the Basis Property of Systems of Root Functions of Regular
	Boundary Value Problems for the Sturm-Liouville Operator. Differ. Equations,
	42, No. 12, 1717-1728 (2006); translation from Differ. Uravn., 42, No. 12,
	1646-1656 (2006) (Russian).

[39] V.A. Il'in. On the unconditional basis property on a closed interval
of systems of eigenfunctions and associated functions of a second-order differential
operator. Dokl. Akad. Nauk {\bf 273} (1983), 1048-1053 (in Russian).

[40] A.S. Makin.  Inverse Problems of Spectral Analysis for the
	Sturm-Liouville Operator with Regular Boundary Conditions: I.
	Differ. Equations, 43, No. 10, 1364-1375 (2007); translation from
	Differ.	 Uravn., 43, No. 10, 1334-1345 (2007) (Russian).

[41] A.S. Makin.  Inverse Problems of Spectral Analysis for the
	Sturm-Liouville Operator with Regular Boundary Conditions: II.
	Differ. Equations, 43, No. 12, 1668-1678 (2007); translation from
	Differ. Uravn., 43, No. 12, 1626-1636 (2007) (Russian).

[42] V.A. Tkachenko. Spectral analysis of non-selfadjoint Hill operator.
 Dokl. Akad. Nauk {\bf 322} (1992), 248-252 (in Russian).

[43] A.S. Makin. A class of boundary value problems for the Sturm-Liouville
operator. Differ. Equations 35, No. 8, 1067-1076 (1999); translation from
 Differ. Uravn. 35, No. 8, 1058-1066 (1999) (Russian).

[44] S. Hoffman. Second-Order Linear Differential Operators Defined by Irregular
Boundary Conditions. Dissertation, Yale University, 1957.

[45] A.S. Makin.  Characterization of the spectrum of the  operator with irregular  boundary conditions. Differ. Equations, No. 10, 1427-1437 (2010);
translation from Differ Uravn., 46, No. 10, 1421-1432 (Russian) (2010).

[46] M.M. Malamud. On completeness on the system of Root Vector Sturm-Liouville
Operator Subject to General Boundary Conditions. Dokl. Akad. Nauk {\bf 419} (2008), 19-22
(in Russian).

[47] M.M. Malamud. On completeness on the system of Root Vector Sturm-Liouville
Operator Subject to General Boundary Conditions. Funct. Anal. Appl. {\bf 42} (2008), pp.198-204.

[48] B.N. Biyarov. On spectral properties of correct restrictions and extensions of the Sturm-Liouville
operator, Differ. Uravn., 30, No. 12, 2027-2032 (Russian) (1994).

[49] A.S. Makin. On the Completeness of the System of root Functions of the Sturm-Liouville Operator with Degenerate Boundary Conditions . Differ. Equations, 50, No. 6, 835-839 (2014); translation from Differ. Uravn., 50, No. 6, 835-838 (Russian) (2014).

[50] A. Minkin. Resolvent growth and Birkhoff-regularity. J. Math. Anal. Appl.
{\bf 323} (2006), 387-402.

[51]  A.S. Makin.  On an Inverse Problem for the Sturm-Liouville Operator with Degenerate Boundary Conditions . Differ. Equations, 50, No. 10, 1402-1406 (2014); translation from Differ. Uravn., 50, No. 10, 1408-1411 (Russian) (2014).

[52]  A.S. Makin. On a two-point boundary value problem for the Sturm-Liouville Operator with nonclassical
spectral asymptotics. Differ. Equations, 49, No. 5, 536-544 (2013); translation from Differ. Uravn., 49, No. 5,
564-572	 (Russian) (2013).

[53]  A.S. Makin.  Problem with Nonclassical Eigenvalue Asymptotics. Differ. Equations, 51, No. 3, 318-324 (2015); translation from Differ. Uravn., 51, No. 3, 317-322 (Russian) (2015).

[54] A.S. Makin. On a new class of boundary value problems for the Sturm-Liouville
operator, Differ. Equations, V. 49, No. 2, 262-266 (2013); translation from Differ. Uravn., 49, No. 2, 260-264 (Russian) (2013).

[55]  A.S. Makin.  On  two-point boundary value problems for the the Sturm-Liouville operator. International Conference "Modern methods of theory of boundary value problems". Voronezh, Russia. May 3-9, 2015. Abstracts of talks, p. 138-140.

[56]  Alexander Makin. Spectral Analysis for the Sturm-Liouville Operator with Degenerate Boundary Conditions. International Conference on Differential and Difference Equations and Applications. Amadora, Portugal, May 18-22, 2015. Abstracts of talks, p. 93-94.

[57]  A.S. Makin. On spectral expansions for the the Sturm-Liouville operator with two-point boundary conditions. International Conference "Function Spaces and Function Approximation Theory". Moscow, Russia. May 25-29, 2015. Abstracts of talks, p. 177.

\medskip
\medskip
\medskip
\medskip

E-mail: alexmakin@yandex.ru

\end{document}